# THE STEPPING STONE MODEL. II: GENEALOGIES AND THE INFINITE SITES MODEL

By Iljana Zähle,[1] J. Theodore Cox[2] and Richard Durrett[3]

*University of Erlangen, Syracuse University and Cornell University*

This paper extends earlier work by Cox and Durrett, who studied the coalescence times for two lineages in the stepping stone model on the two-dimensional torus. We show that the genealogy of a sample of size $n$ is given by a time change of Kingman's coalescent. With DNA sequence data in mind, we investigate mutation patterns under the infinite sites model, which assumes that each mutation occurs at a new site. Our results suggest that the spatial structure of the human population contributes to the haplotype structure and a slower than expected decay of genetic correlation with distance revealed by recent studies of the human genome.

**1. Introduction.** Sequencing of the human genome revealed [see Reich et al. (2001)] a slower decay of linkage disequilibrium (correlation) with distance along chromosomes than predicted by earlier theoretical studies [Kruglyak (1999)]. This correlation is visible in samples as "haplotype structure": sequences can be divided into blocks where there are only a small number of overall mutation patterns (haplotypes); see, for example, Patil et al. (2001). The mapping of genes that cause disease is often done by whole genome association studies that look for regions where there is a correlation between the states of genetic markers and the presence of disease, so it is important to understand the causes of linkage disequilibrium. For surveys, see Ardlie, Kruglyak and Seielstad (2002), Nordborg and Tavaré (2002), and Pritchard and Przeworski (2001). Fixation of beneficial mutations in a population can create haplotype structure [see, e.g., Sabeti et al. (2002)].

Received April 2003; revised October 2003.
[1]Supported by a DAAD fellowship at Cornell University.
[2]Supported in part by NSF Grant 0204422 from the probability program.
[3]Supported in part by NSF Grants 0202935 from the probability program and 0201037 from a joint DMS/NIGMS initiative to support research in mathematical biology.

However, the use of haplotypes from a chromosome 21 region to distinguish multiple prehistoric human migrations [see Jin et al. (1999)] indicates that the spatial structure of the human population plays a role as well.

In this paper we investigate properties of DNA sequences sampled from a population that evolves according to the stepping stone model. Following Cox and Durrett (2002), we represent space as the torus $\Lambda(L)$, which consists of the points in $(-L/2, L/2]^2$ with integer coordinates, and we suppose that at each point $x \in \Lambda(L)$ there is a colony consisting of $N$ diploid or $2N$ haploid individuals, labeled $1, \ldots, 2N$. In contrast to the previous work, we suppose that the population evolves in continuous time, that is, we use the Moran model rather than the one of Wright and Fisher. In a colony with $N$ diploid individuals, the $2N$ copies of the genetic locus are grouped into pairs that are replaced simultaneously. This little bit of realism does not change the properties of the model very much, but adds annoying complications to the proofs, so we follow the common practice of assuming that individuals are a random union of gametes, that is, we suppose our colonies consist of $2N$ haploid individuals.

Ignoring mutations for the moment, in the Moran model each of the individuals in the system is replaced at rate 1. With probability $1-\nu$ ($\nu \in (0,1]$) it is replaced by a copy of an individual that is chosen at random from the colony in which it resides. For convenience we allow the departing individual to be chosen. With probability $\nu$ the departing individual from colony $x$ is replaced by one chosen at random from a nearby colony $y \neq x$ with probability $q(y-x)$, where the difference $y - x \in \Lambda(L)$ is computed componentwise and modulo $L$. Let

$$p(x,y) = (1-\nu)I(x,y) + \nu q(y-x),$$

where $I(x,y) = 1$ if $x = y$ and 0 otherwise. We have separated the kernel into two parts since we are interested in limits as $L \to \infty$ in which the migration rate $\nu$ may converge to 0, but $q(z)$ is a fixed displacement kernel. We suppose $q(z)$ is an irreducible probability distribution on $\mathbb{Z}^2$ with $q((0,0)) = 0$ that has the following properties.

1. $\mathbb{Z}^2$ *symmetry*: $q((x_1, x_2)) = q((-x_1, -x_2))$; $q((x_1, x_2)) = q((x_2, x_1))$.
2. *Finite range*: $q((x_1, x_2)) = 0$ if $\sup_i |x_i| \geq K$ for some $K < \infty$.

We suppose that $L \geq 2K$ so that we do not get confused when we try to define the corresponding random walk transition probability on the torus. The first assumption implies that a single step taken according to $q$ has zero mean and covariance $\sigma^2 I$, where $\sigma^2 = \sum_{x \in \mathbb{Z}^2} x_1^2 q(x) = \sum_{x \in \mathbb{Z}^2} x_2^2 q(x)$. The finite range condition implies $\sigma^2 < \infty$.

To study the behavior of the stepping stone model, we work backwards in time to define a coalescing random walk. When an individual is replaced, its lineage jumps to the one it was replaced by. The history of one individual is thus a random walk. When two lineages come together in one individual they never again separate, so the collection of lineages is a coalescing random walk. As we work backward, let $T_0$ be the amount of time required until the two lineages first reside in the same colony and let $t_0$ be the total amount of time needed for the two lineages to coalesce to one. We begin by considering a sample of size 2, one chosen at random from the colony at 0 and the other an independent choice from the colony at $x$. Let $P_x$ denote the distribution of the genealogy in this case.

Our first result extends Theorem 5 of Cox and Durrett (2002) by giving more refined information about small times. For $0 < \delta \leq 1$ and $c > 0$, let $\Gamma(L, c, \delta) = (L^\delta / \log L, c\delta L^\delta \log L)$.

THEOREM 1. *Suppose that $2N\nu\pi\sigma^2/\log L \to \alpha \in [0, \infty)$ as $L \to \infty$, where $N$ and $\nu$ depend on $L$. For any fixed $\beta_0 > 0$, as $L \to \infty$,*

$$\sup_{\beta_0 \leq \beta \leq \gamma \leq 1} \sup_{|x| \in \Gamma(L,c,\beta)} \left| P_x\left(t_0 > \frac{L^{2\gamma}}{2\nu}\right) - \frac{\beta + \alpha}{\gamma + \alpha} \right| \to 0.$$

If the number of haploid individuals per colony $2N = 1$, $\nu = 1$, and $q$ assigns probability $1/4$ to the four nearest neighbors, then $\alpha = 0$, which is closely related to a result of Cox and Griffeath (1986) for the voter model on $\mathbb{Z}^2$. Indeed their result extends easily to the torus since at time $L^{2\gamma}/2\nu$ with $\gamma < 1$ the particles do not realize they are not on $\mathbb{Z}^2$.

Let $h_L = (1 + \alpha)L^2 \log L/(2\pi\sigma^2\nu)$ and suppose $|x_L| \in \Gamma(L, c, \beta)$. The behavior for larger times as given by Theorem 5 of Cox and Durrett (2002) is

(1.1) $$P_x\left(t_0 > \frac{L^2}{2\nu} + h_L t\right) \to \frac{\beta + \alpha}{1 + \alpha} e^{-t}.$$

Here we have added the term $L^2/2\nu = o(h_L)$ to the Cox and Durrett result so that the times covered by the two results are disjoint. Note that there is a correction to Theorem 5 of Cox and Durrett (2002): In the assumption, $\lim_{L\to\infty} 2N\pi\sigma^2\nu/\log L = \alpha$ has to be replaced by $\lim_{L\to\infty} 4N\pi\sigma^2\nu/\log L = \alpha$. However, in the continuous time model, the first assumption is the correct one.

Our first step in studying the genealogies is to suppose that the random sample is spread out across the torus. Let $\mathcal{G}(L, n, 1)$ be the set of all $n$-point sets where the distance between all points is at least $L/\log L$, that is,

(1.2) $$\mathcal{G}(L, n, 1) = \{A = \{x_1, \ldots x_n\} : \forall i, x_i \in \Lambda(L), \\ \forall i \neq j, |x_i - x_j| \geq L/\log L\}.$$

Let $\zeta_s(A)$ be the coalescing random walk with $\zeta_0 = A$ and let $D_t$ be the pure death process that makes transitions from $k \to k - 1$ at rate $\binom{k}{2}$ with $D_0 = n$. In words, $D_t$ gives the number of lineages at time $t$ in Kingman's coalescent.

THEOREM 2. *Suppose that $2N\nu\pi\sigma^2/\log L \to \alpha \in [0,\infty)$ as $L \to \infty$, where $N$ and $\nu$ depend on $L$. As $L \to \infty$,*

$$\sup_{t\geq 0} \sup_{A\in\mathcal{G}(L,n,1)} |P_A(|\zeta_{h_L t}| = k) - P_n(D_t = k)| \to 0.$$

In the nearest neighbor case with $2N = 1$ this is due to Cox (1989). To express the conclusion in biological terms, we note that in a homogeneously mixing population that consists of a total of $\mathcal{N}$ diploid or $2\mathcal{N}$ haploid individuals, the genealogy on time scale $2\mathcal{N}t$ converges to Kingman's coalescent. Thus for samples with one individual taken from a collection of colonies $A \in \mathcal{G}(L, n, 1)$, our spatial model behaves like a homogeneously mixing population with "effective population size"

$$(1.3) \qquad \mathcal{N}_e = (1+\alpha)\frac{L^2 \log L}{4\pi\sigma^2\nu} \approx NL^2 \cdot \frac{1+\alpha}{2\alpha}.$$

In many genetic studies, sampled individuals are not chosen randomly across the planet. For example, one of the samples in Sabeti et al. (2002) consists of 73 Beni individuals who are civil servants in Benin City, Nigeria. For such a sample, the setup of Theorem 1 is more appropriate. Let $\mathcal{G}(L, n, c, \delta)$ be the set of all $n$-point sets where the distance between all points is in $\Gamma(L, c, \delta)$, that is,

$$(1.4) \quad \mathcal{G}(L,n,c,\delta) = \{A = \{x_1, \ldots x_n\} : \forall i, x_i \in \Lambda(L), \\ \forall i \neq j, |x_i - x_j| \in \Gamma(L,c,\delta)\}.$$

THEOREM 3. *Suppose that $2N\nu\pi\sigma^2/\log L \to \alpha \in [0,\infty)$ as $L \to \infty$, where $N$ and $\nu$ depend on $L$. For any fixed $\beta_0 > 0$, as $L \to \infty$,*

$$\sup_{\beta_0 \leq \beta \leq \gamma \leq 1} \sup_{A \in \mathcal{G}(L,n,c,\beta)} |P_A(|\zeta_{L^{2\gamma}/(2\nu)}| = k) - P_n(D_{\log((\gamma+\alpha)/(\beta+\alpha))} = k)| \to 0,$$

$$\sup_{t\geq 0} \sup_{A \in \mathcal{G}(L,n,c,\beta)} |P_A(|\zeta_{L^2/(2\nu)+h_L t}| = k) - P_n(D_{\log((1+\alpha)/(\beta+\alpha))+t} = k)| \to 0.$$

Again, in the nearest neighbor case with $2N = 1$ the first part is essentially due to Cox and Griffeath (1986). Our result shows that until time $L^2/2\nu$, the particles behave as if they are on $\mathbb{Z}^2$ and then they evolve as predicted by Theorem 2. To prove this result, it is enough to prove the first conclusion and that the configuration at time $L^2/2\nu$ satisfies the assumptions of Theorem 2. The proofs of Theorems 2 and 3 show that when there are $k$ lineages remaining, all $\binom{k}{2}$ pairs have an equal chance to be the next to coalesce, so the partition structure induced by coalescence is the same as in the homogeneously mixing case.

In Section 2 we use Theorems 1–3 to compute various quantities of interest in genetics. Our aim there is to argue that in a population that follows the stepping stone model: (1) genetic correlation decays more slowly with distance along a chromosome than in a homogeneously mixing population and (2) the unusual time scaling before $L^2/2\nu$ can cause haplotype structure. The remainder of the paper is devoted to proofs. Theorems 1–3 are proved in Sections 3–5, respectively.

**2. Applications.** In this section we investigate the impact of spatial structure on the DNA of a sample of $n$ individuals. Since any two humans differ in about 1/1000 nucleotides, we use the infinite sites model which assumes that each new mutation changes a different nucleotide. Some of the formulas we derive are somewhat complicated, so it proves useful to have a concrete example to which to apply our results. The following scenario is motivated by thinking about the human population before it emerged from Africa 100,000 years ago. Our purpose here is not to fit the model to existing data; it is only to show that the stepping stone model can produce patterns that are qualitatively similar to those found in the human genome.

*Concrete example.* Let $L = 100$ and $N = 5$, so the total population size $NL^2 = 50{,}000$. We choose a migration rate $\nu = 0.2$, which corresponds to an average of $N\nu = 1$ migrant per generation, and set $\sigma^2 = 2$. In this case,

$$\alpha \approx \frac{2(5)\pi(2.0)(0.2)}{4.6052} \approx 2.7288,$$

so the effective population size is $(1+\alpha)/2\alpha = 0.68323$ times $NL^2$ or 34,162. To pick a value of $\beta$, we recall Sabeti et al.'s (2002) sample of civil servants in Benin city and somewhat arbitrarily set $\beta = 0.4$.

Theorem 3 implies that if we have a random sample $A \in \mathcal{G}(L, n, c, \beta)$ and change variables

(2.1)
$$\frac{L^{2\gamma}}{2\nu} \to \log\left(\frac{\alpha+\gamma}{\alpha+\beta}\right) \qquad \text{for } \beta \leq \gamma \leq 1,$$
$$\frac{L^2}{2\nu} + (1+\alpha)\frac{L^2 \log L}{2\pi\sigma^2\nu}s \to \log\left(\frac{\alpha+1}{\alpha+\beta}\right) + s \qquad \text{for } s \geq 0,$$

then the genealogy of our sample is that of the ordinary coalescent.

In the example the probability that two lineages do not coalesce by time $L^2/2\nu$ is

$$\frac{\alpha+\beta}{\alpha+1} \approx 0.83909 = 1 - 0.16091,$$

which corresponds to time $\log(1/0.83909) \approx 0.17544$ in the coalescent. If we look at Table 1 in Sabeti et al. (2002), then we see that their sample of 60 Benis produced seven core haplotypes that gave an allelic partition of 14, 13, 10, 10, 9, 3, 1. To compare with our model note that (1) the fraction of pairs that have coalesced is

$$\frac{14(13) + 13(12) + 10(9) + 10(9) + 9(8) + 3(2)}{60 \cdot 59} = \frac{596}{3540} = 0.168$$

and (2) the expected time for a sample of size 60 to be reduced to seven lineages is

$$\sum_{k=8}^{60} \binom{k}{2}^{-1} = \sum_{k=8}^{60} \frac{2}{k-1} - \frac{2}{k} = \frac{2}{7} - \frac{2}{60} = 0.25238.$$

It is useful to reexpress the time change (2.1) in terms of $t/2\nu$ as

$$(2.2) \quad \begin{aligned} \frac{t}{2\nu} &= \frac{L^{2\gamma}}{2\nu} \quad \text{for } \beta \leq \gamma \leq 1 \text{ implies } \gamma = \frac{\log t}{2\log L}, \\ \frac{t}{2\nu} &= \frac{L^2}{2\nu} + (1+\alpha)\frac{L^2 \log L}{2\pi\sigma^2\nu}s \quad \text{implies } s = (t - L^2)\frac{\pi\sigma^2}{(1+\alpha)L^2 \log L}. \end{aligned}$$

Thus for $L^{2\beta} \leq t \leq L^2$,

$$(2.3) \qquad P\left(t_0 \geq \frac{t}{2\nu}\right) \approx \frac{\alpha + \beta}{\alpha + \log t/(2\log L)},$$

while for $t \geq L^2$,

$$(2.4) \qquad P\left(t_0 \geq \frac{t}{2\nu}\right) \approx \frac{\alpha + \beta}{\alpha + 1} \exp\left(-(t - L^2)\frac{\pi\sigma^2}{(1+\alpha)L^2 \log L}\right).$$

*Recombination.* The results above apply to tracing the history of a single nucleotide. To study the decay of genetic correlation with distance we need to investigate the relationship between the genetic history of two different nucleotides separated by a certain distance on a chromosome. To build a mental picture of the process, think of the copies of the first nucleotide as red balls and of the second nucleotide as blue balls. Initially we have $n$ red–blue pairs that represent the initial sample. If we trace back the lineages of the blue balls, then we get a coalescing random walk in which a lineage jumps from $x$ to $y$ when the individual at $x$ is replaced by an offspring of the one at $y$. The same is true for the red balls, but the genealogies of the two colors are coupled. On a given jump, for a red–blue pair, both will be inherited from a single parent with probability $1 - r$ or, with probability $r$, a recombination will occur and the two will be inherited from independently chosen parents. Our next result gives the probability of no recombination before coalescence (NRBC) in a sample of size 2. Let $\ell(u) = \beta \vee \frac{\log(1/u)}{2\log L} \wedge 1$, where $a \vee b = \max\{a,b\}$ and $a \wedge b = \min\{a,b\}$:

$$
\begin{aligned}
P(\text{NRBC}) \approx{} & e^{-rL^{2\beta}/\nu} - e^{-rL^{2\beta}/\nu} - e^{-rL^2/\nu} \frac{\alpha + \beta}{\alpha + \ell(u)} \\
& - e^{-rL^2/\nu} \left(\frac{\alpha + \beta}{\alpha + 1}\right) r \Big/ \left(r + \frac{\pi\sigma^2\nu}{(1+\alpha)L^2 \log L}\right).
\end{aligned}
\tag{2.5}
$$

PROOF OF (2.5). If we condition on $t_0$, $P(\text{NRBC}|t_0) = \exp(-r(2t_0)) = \exp(-(r/\nu)(2\nu t_0))$. Letting $u = r/\nu$ we have

$$
E \exp(-u(2\nu t_0)) = \int_0^\infty e^{-ut} P(2\nu t_0 = t) \, dt. \tag{2.6}
$$

Integrating the above by parts equals

$$
1 - \int_0^\infty u e^{-ut} P(2\nu t_0 \geq t) \, dt.
$$

Using (2.3) and (2.4) and changing variables $t = s + L^2$ in the second integral, we have

$$
\begin{aligned}
& \int_0^\infty u e^{-ut} P(\nu t_0 \geq t) \, dt \\
& \approx 1 - \exp(-uL^{2\beta}) + \int_{L^{2\beta}}^{L^2} u e^{-ut} \frac{\alpha + \beta}{\alpha + \frac{\log t}{2 \log L}} \, dt \\
& \quad + \frac{\alpha + \beta}{\alpha + 1} \int_0^\infty u \exp(-u(L^2 + s)) \exp\left(-s \frac{\pi\sigma^2}{(1+\alpha)L^2 \log L}\right) ds.
\end{aligned}
\tag{2.7}
$$

The last integral is easy to evaluate exactly:

$$
\exp(-uL^2) \left(\frac{\alpha + \beta}{\alpha + 1}\right) u \Big/ \left(u + \frac{\pi\sigma^2}{(1+\alpha)L^2 \log L}\right).
$$

The first integral is

$$
\approx (\exp(-uL^{2\beta}) - \exp(-uL^2))
$$

$$
\times \begin{cases} (\alpha + \beta)/(\alpha + 1), & \text{when } uL^2 \to 0, \\ (\alpha + \beta)/(\alpha + \beta), & \text{when } uL^{2\beta} \to \infty, \\ (\alpha + \beta) \Big/ \left(\alpha + \frac{\log(1/u)}{2 \log L}\right), & \text{otherwise.} \end{cases}
$$

Recalling the definition of $\ell(u)$, and combining this with (2.6) and (2.7), we have

$$
\begin{aligned}
P(\text{NRBC}) \approx{} & \exp(-uL^{2\beta}) - (\exp(-uL^{2\beta}) - \exp(-uL^2)) \frac{\alpha + \beta}{\alpha + \ell(u)} \\
& - \exp(-uL^2) \left(\frac{\alpha + \beta}{\alpha + 1}\right) u \Big/ \left(u + \frac{\pi\sigma^2}{(1+\alpha)L^2 \log L}\right).
\end{aligned}
$$

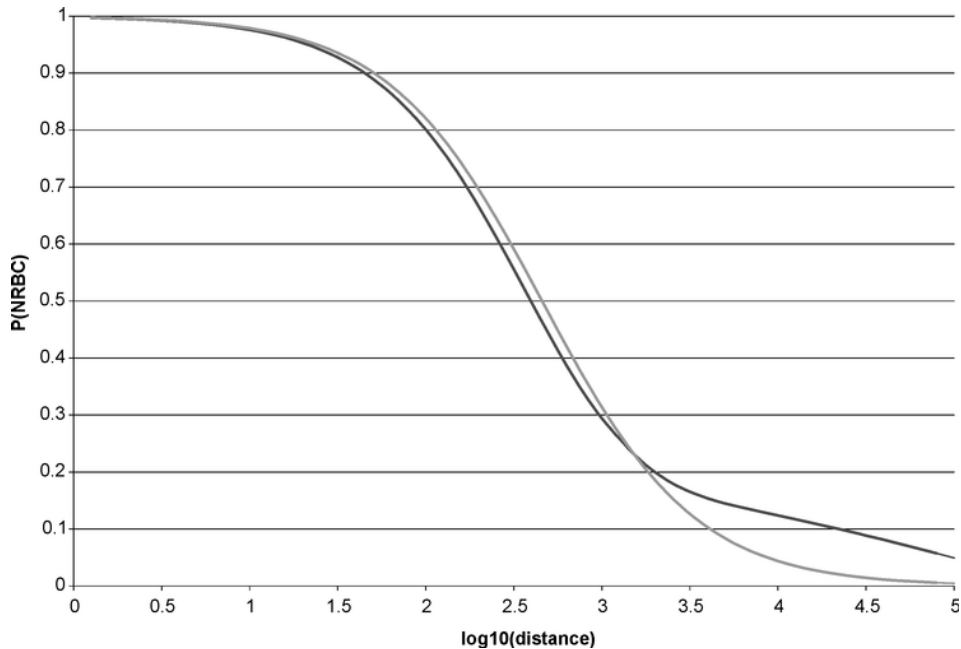

Fig. 1. *Decay of the probability of no recombination before coalescence with the base* 10 *logarithm of distance in the spatial model and in a homogeneously mixing population of size* $\mathcal{N}_e$. *The two are close up to* 1000 *nucleotides, but then the spatial model is larger.*

Since $u = r/\nu$, we have the desired result. □

In our concrete example, $L = 100$ and $\nu = 0.2$, so $qL^2 = 50{,}000r$. Taking $\rho = 10^{-8}$ per nucleotide per generation as a typical value of the recombination rate, we see that the changeover between the second and third terms occurs when the recombination probability between the two nucleotides is $r = 2 \times 10^{-5}$, which corresponds to a distance of 2000 nucleotides. At the other extreme, when $r/\nu = L^{-2\beta}$, the right-hand side is very close to 0. In our example, $\beta = 0.4$ so this occurs for $r = 0.2/100^{-0.8} = 0.0050$, which corresponds to 500,000 nucleotides. Figure 1 shows $P(\text{NRBC})$ for our example for distances 316–100,000 nucleotides and compares it with the result for a homogeneously mixing population of size $\mathcal{N}_e$, defined in (1.3). Note that $P(\text{NRBC})$ is much larger in the spatial model than in the homogeneously mixing case.

*Linkage disequilibrium.* Consider one locus with alleles $A$ and $a$ and a second with alleles $B$ and $b$. A commonly used measure of linkage disequilibrium which is familiar to probabilists is the square of the correlation

coefficient

$$r_{AB}^2 = \frac{(f_{AB} - f_A f_B)^2}{f_A f_a f_B f_b},$$

where $f_c$ is the frequency of genotype $c$. When allele frequencies are larger than 10%, Ohta and Kimura (1971) showed that

$$Er_{AB}^2 \approx \frac{E(f_{AB} - f_A f_B)^2}{E(f_A f_a f_B f_b)} \equiv \sigma_d^2.$$

In a recent paper, McVean (2002) showed that, in general,

$$\sigma_d^2 = \frac{\rho_{ij,ij} - 2\rho_{ij,ik} + \rho_{ij,kl}}{E(T^2)/\text{var}(T) + \rho_{ij,kl}},$$

where $T$ is the coalescence time of a sample of size 2 at one of the loci and the $\rho$'s are correlations between various coalescence times. For example, $\rho_{ij,ik}$ is the correlation of the coalescence time for lineages $i$ and $j$ at locus $x$ with that of lineages $i$ and $k$ at locus $y$, and $i, j, k$ are assumed distinct. For a homogeneously mixing population one can compute [see (12) in McVean (2002)] that

$$\sigma_d^2 = \frac{10 + \rho}{22 + 13\rho + \rho^2}.$$

This calculation [see also Section 2.1 of Durrett (2002)] depends heavily on the fact that the coalescence rates remain constant in time, so we have not been able to calculate this quantity for the stepping stone model. Pritchard and Przeworski (2001) gave simulation results for $r^2$ in a homogeneously mixing population and for population scenarios such as exponential growth and the island model of populations subdivision.

A second commonly used measure of linkage disequilibrium is $D'$, which is defined to be the covariance divided by its maximum possible value. If we suppose without loss of generality that $f_A \geq f_B \geq 1/2$, then

$$D' = \frac{(f_{AB} - f_A f_B)}{f_B - f_A f_B},$$

since in this case the numerator is maximized when $f_{aB} = 0$. Data in Reich et al. (2001) show that $D'$ decays roughly linearly in the logarithm of distance for distances between 5000 and 160,000 nucleotides.

Dawson et al. (2002) studied the decay of $D'$ and $r^2$ with distance for data on human chromosome 21. Their Figure 1 gives results for 1504 markers in which the minor allele frequencies were all greater than 0.2. As the lower two panels show, the average values of $D'$ and $r^2$ do not decay to their limiting values (0 in the case of $r^2$ and 0.2 in the case of $D'$) until the distance is about 200,000 nucleotides. In contrast the upper two panels show that the actual values of $D'$ and $r^2$ for a given pair of markers fluctuate wildly since the values of these statistics depend heavily on where the mutations occur on the genealogical trees. For a more detailed explanation, see Nordborg and Tavaré (2002). Since $D'$ and $r^2$ depend on both the shape of the genealogical tree and the placement of the mutations on it, proving results about these quantities seems difficult.

*Pairwise differences.* If we two individuals at random from a box with side length $L^\beta$, then the average number of places where their DNA sequences differ is $E(2\mu t_0)$, where $\mu$ is the mutation rate for the region under consideration and $t_0$ is the coalescence time of the two lineages. We see below that

$$(2.8) \quad E(2\mu t_0) \approx \frac{\mu}{\nu}\left(\frac{\alpha+\beta}{\alpha+1}\right)\left\{\frac{(1+\alpha)L^2\log L}{\pi\sigma^2} + L^2\left(1 - \frac{1}{2(\alpha+1)\log L}\right)\right\}.$$

Note that the dominant contribution comes from times after $L^2/2\nu$, but, ignoring constants, each successive term is smaller by a factor $1/(\log L) = 0.217$. In our example,

$$E(2\mu t_0) \approx \mu(0.83901)(136{,}646 + 48{,}544) = 155{,}391\mu.$$

Assuming that $\mu = 10^{-8}$, this is $1.55 \times 10^{-3}$, which is in reasonable agreement with the rule of thumb which says that roughly $1/1000$ nucleotides differ between two humans.

PROOF OF (2.8). Using (2.3) and (2.4) in $E(2\mu t_0) = \int_0^\infty P(2\mu t_0 \geq t)\,dt$ we have

$$E(2\mu t_0) \approx \frac{\mu}{\nu}\left\{L^{2\beta} + \int_{L^{2\beta}}^{L^2} \frac{\alpha+\beta}{\alpha + \log t/2\log L}\,dt \right.$$
$$\left. + \left(\frac{\alpha+\beta}{\alpha+1}\right)\int_0^\infty \exp\left(-s\frac{\pi\sigma^2}{(1+\alpha)L^2\log L}\right)ds\right\}.$$

As in the recombination calculation, the second integral is easy to evaluate exactly. To approximate the first we can observe that it is at least $(L^2 - L^{2\beta})(\alpha+\beta)/(\alpha+1)$. For a bound in the other direction we change variables $t = rL^2$ to get

$$\left(\frac{\alpha+\beta}{\alpha+1}\right)L^2\int_{1/L^{2-2\beta}}^{1}\frac{\alpha+1}{\alpha+1+\log r/2\log L}\,dr$$
$$= \left(\frac{\alpha+\beta}{\alpha+1}\right)L^2\int_{1/L^{2-2\beta}}^{1}\frac{dr}{1+\log r/(2(\alpha+1)\log L)}.$$

Using $1 \geq (1+x)(1-x)$ now we have that the above is

$$\geq \left(\frac{\alpha+\beta}{\alpha+1}\right) L^2 \int_{1/L^{2-2\beta}}^{1} \left(1 - \frac{\log r}{2(\alpha+1)\log L}\right) dr$$

$$\approx \left(\frac{\alpha+\beta}{\alpha+1}\right) L^2 \left(1 - \frac{1}{2(\alpha+1)\log L}\right),$$

where in the second step we have used the fact that the antiderivative of $\log r$ is $r \log r - r$ and we have ignored the contribution for the lower limit which is of order $L^{-2(1-\beta)}$. By using the second-order approximation $1/(1+x) \approx 1 - x + x^2$ we can see that the error in the lower bound in the previous display is $O(L^2/(\log L)^2)$. Dropping the smaller term $L^{2\beta}$ and combining our formulas gives the desired result. □

*Larger samples.* To understand properties of larger samples we use the time scale on which the genealogy is the ordinary coalescent, but mutations occur at a time-dependent rate. The first step is to compute the mutation rate. Equations (2.1) and (2.2) together imply that when $L^{2\beta} \leq t \leq L^2$, $\quad \dfrac{t}{2\nu} \to \log\left(\dfrac{\alpha + \log t / 2 \log L}{\alpha + \beta}\right);$

when $t \geq L^2$, $\quad \dfrac{t}{2\nu} \to \log\left(\dfrac{\alpha+1}{\alpha+\beta}\right) + (t - L^2) \dfrac{\pi \sigma^2}{(1+\alpha) L^2 \log L}.$

Setting the right-hand side equal to $u$ and solving, we see that if $u$ is the time variable for the coalescent and $u_1 = \log(\frac{\alpha+1}{\alpha+\beta})$, then

when $0 \leq u \leq u_1$, $\quad t = \exp([(\alpha+\beta)e^u - \alpha](2\log L));$

when $u \geq u_1$, $\quad t = L^2 + \left\{u - \log\left(\dfrac{\alpha+1}{\alpha+\beta}\right)\right\} \dfrac{(1+\alpha)L^2 \log L}{\pi \sigma^2}.$

Differentiating we have

when $0 \leq u \leq u_1$, $\quad \dfrac{dt}{du} = t(u)(\alpha+\beta)e^u (2 \log L);$

when $u \geq u_1$, $\quad \dfrac{dt}{du} = \dfrac{(1+\alpha)L^2 \log L}{\pi \sigma^2}.$

In the second time interval the mutation rate is constant and has rate

$$\frac{\mu}{2\nu} \cdot \frac{(1+\alpha)L^2 \log L}{\pi \sigma^2}.$$

The first time interval is the set of $u_\gamma = \log((\alpha+\gamma)/(\alpha+\beta))$ with $\beta \leq \gamma \leq 1$. At these times we have $t(u_\gamma) = L^{2\gamma}$ and hence mutation rate

$$\frac{\mu}{2\nu} \cdot (\alpha+\gamma) L^{2\gamma} \log L.$$

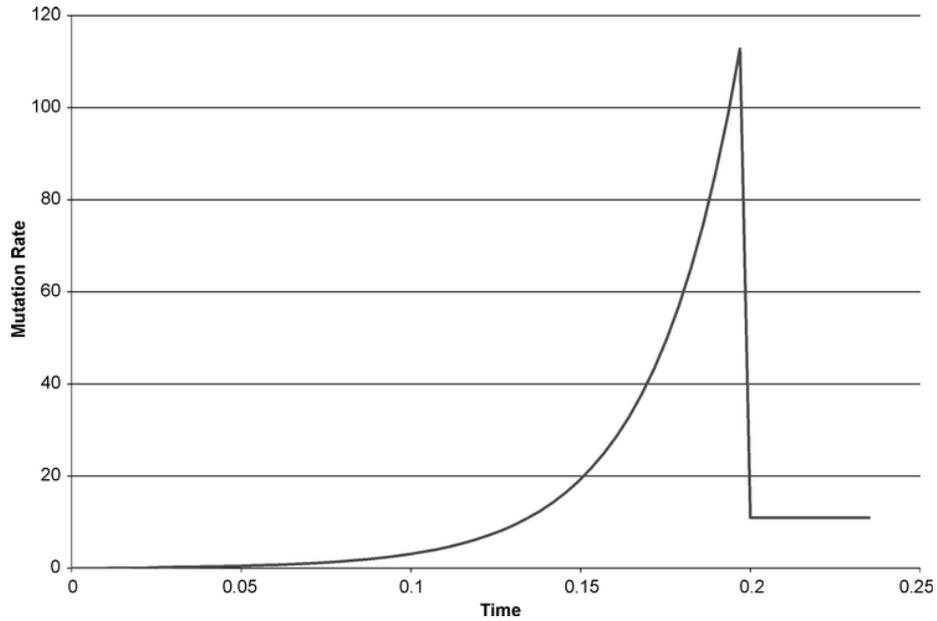

Fig. 2. *Mutation rate in the coalescent as a function of time.*

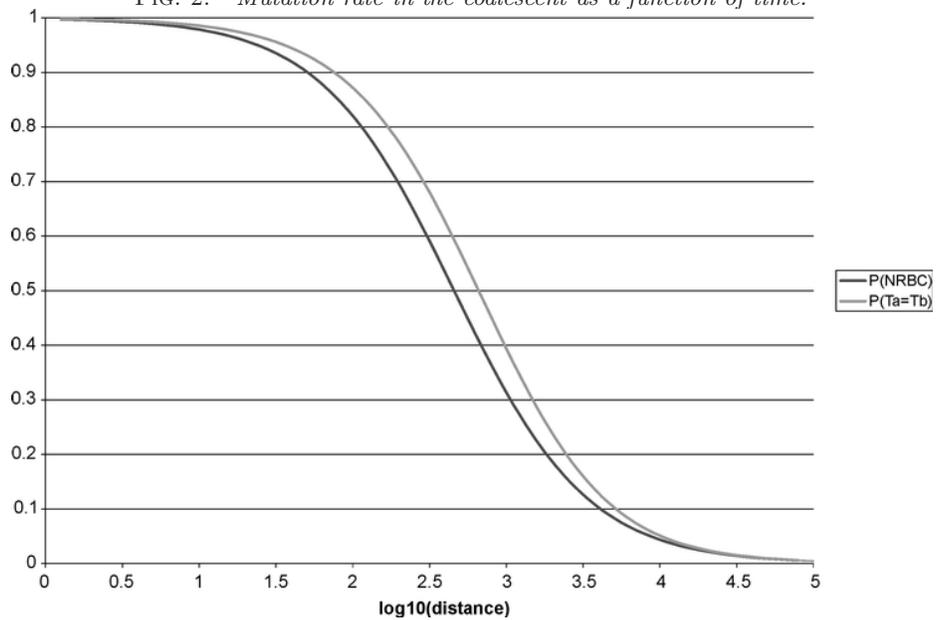

Fig. 3. *Probability of no recombination before coalescence compared to the probability that the coalescence time of the a locus is equal to that of the b locus in the homogeneously mixing case. In our case there should be a more substantial difference since a recombination will put the a and b loci which just separated into the same or nearby colonies.*

To see what this means, suppose that the mutation rate is $\mu = 10^{-8}$ per nucleotide and consider a region with 10,000 nucleotides [roughly the size of the core haplotypes in the G6PD example in Sabeti et al. (2002)]. Then using the calculation after (2.8), the mutation rate is

$$\text{for } u \geq u_1, \qquad 10^{-4} \cdot \frac{136{,}646}{2} = 6.83;$$

when $u = u_\gamma$, $\quad \dfrac{10^{-4}}{0.4} \cdot (2.7287 + \gamma) 10^{4\gamma}(4.6051) = (15.71 + 5.75\gamma)10^{4(\gamma-1)}.$

When $\gamma = 1$ the rate is 21.46. There is a discontinuity in the rate at $u_1$ due to the different ways in which the process is scaled for $t \leq L^2/(2\nu)$ and $t \geq L^2/(2\nu)$. The rate is very large at the end of the first interval, but is large for only a short time. For a picture, see Figure 2. By calculations after (2.8), for a sample of size 2 from a region with 10,000 nucleotides, an average of 4.07 mutations occur before $u_1$ and an average of 11.46 occur after $u_1$. The previous calculation shows that those that occur before $u_1$ occur close to that time. Since the rate decays exponentially fast as we move back toward time 0, this suggests that in a large sample, the first mutations occur after a considerable amount of coalescence has occurred, leading to large sets of individuals with identical mutation patterns (i.e., haplotype structure in the data).

**3. Proof of Theorem 1.** Let $\kappa_L = 1 - (2 \log \log L)/\log L$ and recall that $\Gamma(L, c, \delta) = (L^\delta/\log L, c\delta L^\delta \log L)$. Our first step is to show that up to time $L^{2\kappa_L}/(2\nu) = L^2/(2\nu(\log L)^4)$ the particles do not know that they are on the torus. We then show that (a) if $\beta \leq \kappa_L$, then the probability $t_0$ occurs between times $L^{2\kappa_L}/(2\nu)$ and $L^2/(2\nu)$ is small, and (b) if $\kappa_L \leq \beta \leq 1$, the probability $t_0$ occurs before time $L^2/(2\nu)$ is small.

By rotation invariance we can suppose without loss of generality that $0 \in A$. We suppose that our random walks $X_t$ on the torus are constructed from a random walk $W_t$ on $\mathbb{Z}^d$ (with kernel $p$ and jump rate 1) so that $X_t = W_t \bmod L$. Let $P_x$ denote the probability distribution when the random walk is started in $x$. Note that the variance of $p$ is $\nu\sigma^2$. Using the $L^2$ maximal inequality for martingales, $(a+b)^2 \leq 4(a^2 + b^2)$ and $|x_i| \leq c\beta L^\beta \log L$, and then $\beta \leq \kappa_L$, we can estimate that for $x_i \in A \in \mathcal{G}(L, n, c, \delta)$,

$$\begin{aligned}
P_{x_i}\left(\max_{0 \leq t \leq L^2/(2\nu(\log L)^4)} |W_t| > \frac{L}{3}\right) \\
\leq \frac{C}{L^2} E_{x_i}[|W_{L^2/(2\nu(\log L)^4)}|^2] \\
\leq \frac{C}{L^2}(|x_i|^2 + E_{x_i}[|W_{L^2/(2\nu(\log L)^4)} - x_i|^2]) \\
\leq \frac{C}{L^2}\left(\frac{c^2 L^2 (\log L)^2}{(\log L)^4} + \frac{\sigma^2 L^2}{(\log L)^4}\right),
\end{aligned} \qquad (3.1)$$

which converges to 0 as $L \to \infty$. (Here and in similar estimates below the constant $C$ may change from line to line.) This means we can study the system on $\mathbb{Z}^2$.

We begin with some preliminary results for random walks on $\mathbb{Z}^2$. Many of these facts and their proofs are standard. We give the details because we need to know the results are uniform in various parameters. Let $\bar{X}_t = W_t^1 - W_t^2$ be the difference of two independent continuous time random walks with kernel $p$ and jump rate 1. Since $p$ is symmetric, $\bar{X}_t$ is a continuous time random walk with kernel $p$ and jump rate 2. Taking the special form of $p$ into account, we define $Y_t = \bar{X}_{t/(2\nu)}$, which is a continuous time random walk with kernel $q$ and jump rate 1. Let $\bar{T}_0 = \inf\{t \geq 0 : \bar{X}_t = 0\}$ be the first hitting time of the origin and let $T_0^* = \inf\{t \geq 0 : Y_t = 0\}$ be the corresponding time for $Y$. Since a trivial time change separates the two processes, we can study either one. In general we choose to study $Y_t$, which has the annoying factor $\nu$ eliminated. Recall that $P_0$ denotes the probability distribution when the random walk is started in 0. By $P_q$ we mean that the starting point is chosen according to $q$.

LEMMA 3.1. *As $t \to \infty$,*

$$P_q(T_0^* > t) \sim \frac{2\pi\sigma^2}{\log t}.$$

PROOF. Decomposing according to the last visit to zero before time $t$ (more precisely, the leaving time of the last visit),

$$1 = \int_0^t P_0(Y_s = 0) P_q(T_0^* > t - s) \, ds + P_0(Y_t = 0).$$

Dropping the $-s$ we have

(3.2) $$P_q(T_0^* > t) \leq \frac{1}{\int_0^t P_0(Y_s = 0) \, ds} \sim \frac{2\pi\sigma^2}{\log t}.$$

That last statement can be seen as follows. The local central limit theorem gives

$$\lim_{s \to \infty} (2\pi s P_0(Y_s = 0) - 1/\sigma^2) = 0.$$

Integrating this yields $\int_0^t P_0(Y_s = 0) \, ds \sim 2\pi\sigma^2/\log t$. A continuous time version of the local central limit theorem can be found, for instance, in Zähle [(2002), Proposition D.2].

For the lower bound we decompose by the last visit to zero before time $t + t \log t$ and compute as before; hence

$$1 = \int_0^{t+t\log t} P_0(Y_s = 0) P_q(T_0^* > t + t\log t - s) \, ds + P_0(Y_{t+t\log t} = 0).$$

We split the integral at time $t\log t$. In the first part we estimate

$$P_q(T_0^* > t + t\log t - s) \leq P_q(T_0^* > t)$$

and in the second part we estimate this probability by 1. We end up with

$$(3.3) \quad P_q(T_0^* > t) \geq \frac{1 - \int_{t\log t}^{t+t\log t} P_0(Y_s = 0)\,ds - P_0(Y_{t+t\log t} = 0)}{\int_0^{t\log t} P_0(Y_s = 0)\,ds}.$$

Let $I(s,t) = \int_s^t P_0(Y_r = 0)\,dr$. Again by the local central limit theorem, $I(0, t\log t) \sim \log t/2\pi\sigma^2$, while

$$I(t\log t, t + t\log t) \sim \frac{1}{2\pi\sigma^2}\log\left(1 + \frac{1}{\log t}\right) \to 0.$$

This completes the proof. $\square$

LEMMA 3.2. *Given $\beta_0$, there exists a constant $C$ so that for all $L \geq L_0$ and $\beta_0 \leq \gamma \leq 1$,*

$$P_q\left(\frac{L^{2\gamma}}{2(\log L)^{3/2}} \leq T_0^* < L^{2\gamma}\right) \leq \frac{C\log\log L}{(\log L)^2}.$$

PROOF. Let $u_1 = L^{2\gamma}/2(\log L)^{3/2}$ and $u_2 = L^{2\gamma}$. By (3.2) and (3.3),

$$P_q(u_1 < T_0^* \leq u_2)$$
$$\leq \frac{1}{I(0, u_1)} - \frac{1 - I(u_2 \log u_2, u_2 + u_2 \log u_2)}{I(0, u_2 \log u_2)}$$
$$= \frac{I(u_1, u_2 \log u_2) + I(0, u_1)I(u_2 \log u_2, u_2 + u_2 \log u_2)}{I(0, u_1)I(0, u_2 \log u_2)}.$$

Using the local central limit theorem,

$$I(0, u_1) \sim \frac{\log u_1}{2\pi\sigma^2} \sim \frac{\gamma \log L}{\pi\sigma^2},$$
$$I(0, u_2 \log u_2) \sim \frac{\log(u_2 \log u_2)}{2\pi\sigma^2} \sim \frac{\gamma \log L}{\pi\sigma^2},$$
$$I(u_1, u_2 \log u_2) \sim \frac{\log(u_2 \log u_2) - \log u_1}{2\pi\sigma^2}$$
$$= \frac{\log(2\gamma \log L) + \log(2(\log L)^{3/2})}{2\pi\sigma^2} \sim \frac{5\log\log L}{4\pi\sigma^2},$$
$$I(u_2 \log u_2, u_2 + u_2 \log u_2) \to 0.$$

Plugging these results into the previous formula gives the result. $\square$

Let $R_0 = 0$ and, for $k \geq 1$, let $Q_k$ be the first time the random walk $Y_t$ leaves colony 0 after time $R_{k-1}$ and let $R_k$ be the first hitting time of 0 after time $Q_k$, that is,

$$Q_k = \inf\{s > R_{k-1} : Y_s \neq 0\},$$
$$R_k = \inf\{s > Q_k : Y_s = 0\},$$

and let $K = \min\{k \geq 1 : R_k - Q_k > L^{2\gamma}\}$. Then $K$ is geometric with success probability

$$\vartheta = P_q(T_0^* > L^{2\gamma}).$$

Consider $Y_t$ as a random walk with jump rate $1/\nu$ and jump kernel $p$ and let $N_k$ be the number of jumps that land in colony 0 at times in $[R_{k-1}, Q_k)$. The $N_k$ are independent and are geometric with success probability $\nu$. Define $\mathcal{O}_L$ to be the number of jumps that land in colony 0 before time $L^{2\gamma}$ and let $\mathcal{O}_K = \sum_{k=1}^K N_k$. We are interested primarily in $\mathcal{O}_L$, but $\mathcal{O}_K$ is easier to analyze since it is a sum of independent random variables. The next result shows that $\mathcal{O}_L = \mathcal{O}_K$ with high probability.

LEMMA 3.3. *Given $0 < \beta_0 < 1$ fixed, there is a constant $C$ so that for $\beta_0 \leq \gamma \leq 1$ and $L \geq L_0$,*

$$P_0(\mathcal{O}_K \neq \mathcal{O}_L) \leq C \left( \frac{\log L}{L^{2\beta_0}} + \frac{1}{\sqrt{\log L}} + \frac{\log \log L}{\sqrt{\log L}} \right).$$

PROOF. Since $R_K - Q_K > L^{2\gamma}$, it is enough to bound $P_0(Q_K \geq L^{2\gamma})$. We decompose

$$(3.4) \qquad Q_K = \sum_{k=1}^K (Q_k - R_{k-1}) + \sum_{k=1}^{K-1} (R_k - Q_k).$$

Note that $Q_k - R_{k-1}$ is a sum of $N_k$ independent exponential variables with mean $\nu$. Hence $E_0[Q_k - R_{k-1}] = 1$. For the first sum on the right-hand side of (3.4) we use Markov's inequality and Lemma 3.1 to conclude

$$P_0 \left( \sum_{k=1}^K (Q_k - R_{k-1}) \geq \frac{L^{2\gamma}}{2} \right) \leq \frac{2}{L^{2\gamma}} E_0 \left[ \sum_{k=1}^K (Q_k - R_{k-1}) \right]$$
$$= \frac{2}{\vartheta L^{2\gamma}} \leq C \frac{\gamma \log L}{L^{2\gamma}} \leq C \frac{\log L}{L^{2\beta_0}}.$$

For the second sum in (3.4) note that if $K < (\log L)^{3/2}$ and $G = \{R_k - Q_k < L^{2\gamma}/(2(\log L)^{3/2})$ for all $k < K\}$, then

$$\sum_{k=1}^{K-1} (R_k - Q_k) < \frac{L^{2\gamma}}{2}.$$

Next, by Markov's inequality and by Lemma 3.1,
$$P_0(K \geq (\log L)^{3/2}) \leq \frac{E_0 K}{(\log L)^{3/2}} = \frac{1}{\vartheta(\log L)^{3/2}} \leq \frac{C}{\sqrt{\log L}}.$$

Furthermore, since $R_k - Q_k \leq L^{2\gamma}$ for $k < K$, using Lemma 3.2 gives
$$P_0(G^c \cap \{K < (\log L)^{3/2}\}) \leq (\log L)^{3/2} P_q\left(\frac{L^{2\gamma}}{2(\log L)^{3/2}} < T_0^* \leq L^{2\gamma}\right)$$
$$\leq \frac{C \log \log L}{\sqrt{\log L}}.$$

Combining our estimates gives the indicated result. □

We are now ready to start to estimate the time for two lineages to coalesce. The first step is to consider the coalescence time when they start in the same colony. Then we study the time required to come to the same colony.

LEMMA 3.4. *If $2N\pi\sigma^2\nu/\log L \to \alpha$, then as $L \to \infty$,*
$$\sup_{\beta_0 \leq \gamma \leq \kappa_L} \left| P_0(t_0^* > L^{2\gamma}) - \frac{\alpha}{\alpha + \gamma} \right| \to 0.$$

PROOF. Since the probability of coalescence when two lineages land in the same colony is $1/2N$,
$$P_0(t_0^* > L^{2\gamma}) = E_0\left(1 - \frac{1}{2N}\right)^{\mathcal{O}_L}.$$

Since $\mathcal{O}_L \leq \mathcal{O}_K$, we have
$$0 \leq E_0\left(1 - \frac{1}{2N}\right)^{\mathcal{O}_L} - E_0\left(1 - \frac{1}{2N}\right)^{\mathcal{O}_K} \leq P_0(\mathcal{O}_L \neq \mathcal{O}_K) \to 0$$

by Lemma 3.3. Since $\mathcal{O}_K$ is geometric with success probability $\vartheta\nu$,
$$E_0\left(1 - \frac{1}{2N}\right)^{\mathcal{O}_K} = \sum_{k=1}^{\infty} \vartheta\nu(1 - \vartheta\nu)^{k-1}\left(1 - \frac{1}{2N}\right)^k$$
$$= \frac{\vartheta\nu(1 - 1/2N)}{\vartheta\nu(1 - 1/2N) + (1/2N)}.$$

By Lemma 3.1, $2N\vartheta\nu \to \alpha/\gamma$ uniformly for $\beta_0 \leq \gamma \leq 1$, which completes the proof. □

LEMMA 3.5. *For any fixed $\rho > 0$, there exists a constant $C_\rho$ so that for all $x$ and $u \geq u_0$, where $u_0 < \infty$,*
$$P_x(Y_s = 0 \text{ for some } s \in [u/(\log u)^\rho, u)) \leq \frac{C_\rho \log \log u}{\log u}.$$

PROOF. By considering time $\tau_0$ of the first visit to 0 after time $u/(\log u)^\rho$ we have

$$\int_{u/(\log u)^\rho}^{2u} P_x(Y_s = 0)\,ds = \int_{u/(\log u)^\rho}^{2u} P_x(\tau_0 \in dt) \int_0^{2u-t} ds\, P_0(Y_s = 0).$$

Now we replace the first integral on the right-hand side by $\int_{u/(\log u)^\rho}^{u}$ and then replace the second integral by $\int_0^u$. This yields the estimate

$$P_x(Y_s = 0 \text{ for some } s \in [u/(\log u)^\rho, u)) \cdot \int_0^u P_0(Y_s = 0)\,ds$$

$$\leq \int_{u/(\log u)^\rho}^{2u} P_x(Y_s = 0)\,ds.$$

The local central limit theorem shows that if $\phi$ is the limiting normal density function, then

$$\sup_{x \in \mathbb{Z}^2} |sP_0(Y_s = x) - \phi(x/s^{1/2})| \to 0 \qquad \text{as } s \to \infty.$$

From this it follows that if $u \geq u_0$, the probability of interest is bounded by

$$\frac{C \int_{u/(\log u)^\rho}^{2u} 1/s\,ds}{\log u} \leq \frac{C_\rho \log \log u}{\log u},$$

which gives the desired result. □

LEMMA 3.6. *There exists*

$$\sup_{\beta_0 \leq \beta \leq \gamma \leq \kappa_L} \sup_{|x| \in \Gamma(L,c,\beta)} \left| P_x(T_0^* \leq L^{2\gamma}) - \left(1 - \frac{\beta}{\gamma}\right) \right| \to 0.$$

PROOF. The $L^2$ maximal inequality for martingales implies that

$$P_y(T_0^* \leq |y|^2/\log|y|) \leq P_0\left(\max_{0 \leq t \leq |y|^2/\log|y|} |Y_t| \geq |y|\right) \leq C/(\log|y|).$$

Using this result with Lemma 3.5 for $u = |y|^2(\log|y|)^5$ and $\rho = 6$ it follows that

$$(3.5) \qquad P_y(T_0^* \leq |y|^2(\log|y|)^5) \leq \frac{C}{\log|y|} + \frac{C' \log\log|y|}{\log|y|}.$$

Recalling that $\Gamma(L, c, \beta) = (L^\beta/\log L, c\beta L^\beta \log L)$, we have for $\beta_0 \leq \beta \leq 1$ and $L$ large enough

$$\sup_{|y| \in \Gamma(L,c,\beta)} P_y(T_0^* \leq |y|^2(\log|y|)^5) \leq \frac{C}{\log L} + \frac{C' \log\log L}{\log L}.$$

Repeating the reasoning from the proof of Lemmas 3.1 and 3.5 shows that

$$p^* \equiv P_x(Y_s = 0 \text{ for some } s \in [L^{2\beta}(\log L)^5, L^{2\gamma}])$$

$$\geq \frac{\int_{L^{2\beta}(\log L)^5}^{L^{2\gamma}} P_x(Y_s = 0)\, ds}{\int_0^{L^{2\gamma}} P_0(Y_s = 0)\, ds}.$$

In the other direction,

$$p^* \leq \frac{\int_{L^{2\beta}(\log L)^5}^{2L^{2\gamma}} P_x(Y_s = 0)\, ds}{\int_0^{L^{2\gamma}} P_0(Y_s = 0)\, ds}.$$

The local central limit theorem implies that

$$\sup_{\beta_0 \leq \beta \leq \gamma \leq \kappa_L} \sup_{|x| \in \Gamma(L,c,\beta)} \sup_{s \in [L^{2\beta}(\log L)^5, 2L^{2\gamma}]} |sP_x(Y_s = 0) - 1/2\pi\sigma^2| \to 0.$$

Combining these estimates we have that if $\beta_0 \leq \beta \leq \gamma \leq \kappa_L$ and $\varepsilon > 0$, then if $L$ is large,

$$(1-\varepsilon)\frac{2\gamma \log L - 2\beta \log L - 5\log\log L}{2\gamma \log L}$$

$$\leq p^* \leq (1+\varepsilon)\frac{\log 2 + 2\gamma \log L - 2\beta \log L - 5\log\log L}{2\gamma \log L}$$

and we have the desired result. □

The final step is to combine Lemmas 3.4–3.6.

PROOF OF THEOREM 1. Recall that $T_0^* = \inf\{t \geq 0 : Y_t = 0\}$, where $Y_t$ is a continuous time random walk on $\mathbb{Z}^2$ with kernel $q$ and jump rate 1. That means $T_0^*$ is the time two lineages need to come to the same colony but after a time change with $2\nu$ in the system on $\mathbb{Z}^2$. Let $t_0^*$ be the coalescing time after the same time change in the system on $\mathbb{Z}^2$. Decomposing according to the value of $T_0^*$,

$$P_x(t_0^* > L^{2\gamma}) = P_x(T_0^* > L^{2\gamma}/2, t_0^* > L^{2\gamma}) + P_x(T_0^* \leq L^{2\gamma}/2, t_0^* > L^{2\gamma}).$$

For the first term on the right-hand side we note that if $L^{2\beta_0} \geq u_0$, then Lemma 3.5 with $\rho = 1$ implies that for $\beta_0 \leq \beta \leq \gamma \leq \kappa_L$ and all $x$,

$$0 \leq P_x\left(T_0^* > \frac{L^{2\gamma}}{2}, t_0^* > L^{2\gamma}\right) - P_x(T_0^* > L^{2\gamma})$$

$$\leq P_x\left(\frac{L^{2\gamma}}{2} < T_0^* \leq L^{2\gamma}\right) \leq \frac{C\log\log L}{\log L}.$$

For the second term we note that if $L^{2\beta_0} \geq u_0$, then Lemma 3.5 implies

$$0 \leq P_x\left(T_0^* \leq \frac{L^{2\gamma}}{2}, t_0^* > L^{2\gamma}\right) - P_x\left(T_0^* \leq \frac{L^{2\gamma}}{2}\right) P_0(t_0^* > L^{2\gamma})$$

$$\leq P_x\left(T_0^* \leq \frac{L^{2\gamma}}{2}\right) P_0\left(\frac{L^{2\gamma}}{2} < t_0^* \leq L^{2\gamma}\right) \leq \frac{C \log \log L}{\log L}.$$

Using Lemmas 3.4 and 3.6 now it follows that

$$\sup_{\beta_0 \leq \beta \leq \gamma \leq \kappa_L} \sup_{|x| \in \Gamma(L,c,\beta)} \left| P_x(t_0^* > L^{2\gamma}) - \left(\frac{\beta}{\gamma} + \left(1 - \frac{\beta}{\gamma}\right)\frac{\alpha}{\alpha + \gamma}\right) \right| \to 0,$$

which is the desired result up to $\kappa_L$.

It remains to show that (a) if $\beta \leq \kappa_L = 1 - (2 \log \log L)/\log L$, then the probability $t_0$ occurs between times $L^2/(2\nu(\log L)^4)$ and $L^2/(2\nu)$ is small, and (b) if $\kappa_L \leq \beta \leq 1$, the probability $t_0$ occurs before time $L^2/(2\nu)$ is small. Let $\hat{X}_t = X_t^1 - X_t^2$ be the difference random walk of two independent continuous time random walks on the torus with kernel $p$ and jump rate 1 and let $\hat{Y}_t = \hat{X}_{t/(2\nu)}$. Let $\hat{T}_0$ and $\hat{T}_0^*$ be the hitting times of 0 for $\hat{X}_t$ and $\hat{Y}_t$.

LEMMA 3.7. *There is a constant $C$ so that for all $L$ and $x \in \Lambda(L)$,*

$$P_x(\hat{Y}_s = 0) \leq \frac{C}{s \wedge L^2}.$$

PROOF. This is straightforward given the estimates in the Appendix of Cox and Durrett (2002). First consider $s \leq L^2$. In this case one can use a local central limit theorem from Bhattacharya and Rao (1976) for random walks on $\mathbb{Z}^2$ and sum over $zL^2$ for $z \in \mathbb{Z}^2$ to prove the result. The result extends to $s \geq L^2$ by noting that the Markov property implies that the largest value of $P_x(\hat{Y}_s = 0)$ is decreasing in $s$. □

Using Lemma 3.7 and repeating the proof of Lemma 3.5 shows:

LEMMA 3.8. *If $L \geq L_0$, then*

$$P_x\left(\hat{Y}_s = 0 \text{ for some } s \in \left[\frac{L^2}{(\log L)^4}, L^2\right]\right) \leq \frac{C_\rho \log \log L}{\log L}.$$

PROOF. By considering the first visit to 0 after time $L^2/(\log L)^4$ we have

$$P_x(\hat{Y}_s = 0 \text{ for some } s \in [L^2/(\log L)^4, L^2]) \cdot \int_0^{L^2} P_0(\hat{Y}_s = 0) \, ds$$

$$\leq \int_{L^2/(\log L)^4}^{2L^2} P_x(\hat{Y}_s = 0) \, ds.$$

Lemma 3.7 gives an upper bound on the right-hand side. To get a lower bound on the integral that involves $P_0$, we stop at time $L^2/(\log L)^4$. The estimate in (3.1) shows that up to this time the random walk does not realize it is not on $\mathbb{Z}^2$, so using the local central limit theorem we conclude that if $L \geq L_0$, the probability of interest is bounded by

$$\frac{C \int_{L^2/(\log L)^4}^{2L^2} 1/(s \wedge L^2)\ ds}{\log L} \leq \frac{C_\rho \log \log L}{\log L},$$

which gives the desired result. □

Lemma 3.8 gives (a). To establish (b) now, we note that arguing as in the proof of (3.5) but using Lemma 3.5 with $\rho = 7$ gives

$$P_x(\hat{T}_0^* \leq |x|^2 (\log |x|)^6) \leq \frac{C \log \log L}{\log L}.$$

If $|x| \in \Gamma(L, c, \beta)$ and $\beta \geq \kappa_L$, then $|x| \geq L^{2\kappa_L}/\log L \geq L^2/(\log L)^5$, so if $L \geq L_0$, it follows that $|x|^2(\log |x|)^6 \geq L^2$. This establishes (b) and the proof of Theorem 1 is complete. □

**4. Proof of Theorem 2.** Recall that $h_L = (1 + \alpha)L^2 \log L/(2\pi\sigma^2 \nu)$ and

$$\mathcal{G}(L, n, 1) = \{A = \{x_1, \ldots, x_n\} : \forall i, x_i \in \Lambda(L), \forall i \neq j, |x_i - x_j| \geq L/\log L\}.$$

Theorem 5 of Cox and Durrett (2002) gives the asymptotic behavior of the coalescence time of two particles that are separated by $L/\log L$. The key to deriving a result for the genealogy is to show that when two particles coalesce, the others are separated. Recall that $\zeta_t$ is the system of lineages. Now $\zeta_t$ is started in $A = \{x_1, \ldots, x_4\} \in \mathcal{G}(L, n, 4)$. By $\zeta_t(x_i)$ we denote the position at time $t$ of the lineage started in $x_i$. Let $\tau_{ij}$ be the coalescing time of the two lineages started in $x_i$ and $x_j$ and let $\tau$ be the minimum of $\tau_{ij}$ ($i \neq j$).

LEMMA 4.1. *Let $\zeta$ be started with four lineages in $A = \{x_1, \ldots, x_4\}$. As $L \to \infty$, uniformly in $A \in \mathcal{G}(L, 4, 1)$,*

(4.1) $$\int_0^\infty P\left(\tau = \tau_{12} \in ds, |\zeta_s(x_1) - \zeta_s(x_3)| \leq \frac{L}{\log L}\right) \to 0,$$

(4.2) $$\int_0^\infty P\left(\tau = \tau_{12} \in ds, |\zeta_s(x_3) - \zeta_s(x_4)| \leq \frac{L}{\log L}\right) \to 0.$$

PROOF. The proof is a modification of the proof of (3.5) of Cox (1989). As in his paper, we just prove the first result and leave it to the reader to check that the same proof with small changes gives the second result. Let $(X_t(x_i))_{t\geq 0}$, $i = 1, \ldots, 4$, be independent random walks on $\Lambda(L)$ with kernel $p$ and jump rate 1. Then

$$\int_0^\infty P\left(\tau = \tau_{12} \in ds, |\zeta_s(x_1) - \zeta_s(x_3)| \leq \frac{L}{\log L}\right)$$
$$\leq P(\tau = \tau_{12} \leq t_L h_L)$$
$$+ \int_{t_L h_L}^\infty \sum_{y,z:\, |y-z|\leq L/\log L} P(\tau_{12} \in ds, X_s(x_1) = y) P(X_s(x_3) = z).$$

If $t_L = 1/\log L$, the first term on the right-hand side tends to 0 by Theorem 5 in Cox and Durrett (2002); see (1.1) but remove the added term $L^2/2\nu$. By the estimate in Lemma 3.7, the sum over $z$ in the second term is at most

$$\left(\frac{L}{\log L}\right)^2 \cdot \frac{C}{L^2} \to 0.$$

Since

$$\int_{t_L h_L}^\infty \sum_y P(\tau_{12} \in ds, X_s(x_1) = y) \leq 1,$$

the desired result follows. □

Recall

$$\mathcal{G}(L, n, c, \delta) = \{A = \{x_1, \ldots, x_n\} : \forall i, x_i \in \Lambda(L), \forall i \neq j, |x_i - x_j| \in \Gamma(L, c, \delta)\},$$

where $\Gamma(L, c, \delta) = (L^\delta/\log L, c\delta L^\delta \log L)$.

LEMMA 4.2. *If $2N\nu\pi\sigma^2/\log L \to \alpha \in [0, \infty)$ as $L \to \infty$, where $N$ and $\nu$ depend on $L$, then as $L \to \infty$,*

$$\sup_{t\geq 0} \sup_{A \in \mathcal{G}(L,n,c,2)} \left| P(|\zeta_{h_L t}(A)| = n) - \exp\left(-\binom{n}{2} t\right) \right| \to 0.$$

PROOF. Since the two quantities are monotone decreasing in $t$, it suffices to prove the result for each fixed $t$. The proof is a modification of the proof of (3.1) in Cox (1989). We need the notation, $H_t(i,j) = \{\tau_{ij} \leq h_L t\}$, $F_t(i,j) = \{\tau = \tau_{ij} \leq h_L t\}$ and $q(t) = P(\tau \leq h_L t)$. We decompose $H_t(i,j)$:

(4.3)
$$P(H_t(i,j)) = P(\tau = \tau_{ij} \leq h_L t)$$
$$+ \sum_{\{k,l\}\neq\{i,j\}} \int_0^{h_L t} P(\tau = \tau_{kl} \in ds, \tau_{ij} \leq h_L t).$$

The $k, l$ term in the second sum is

$$= \int_0^{h_L t} \sum_{y,z} P(\tau = \tau_{kl} \in ds, X_s(x_i) = y, X_s(x_j) = z) P(|\zeta_{h_L t - s}(\{y,z\})| = 1).$$

By Lemma 4.1 we can neglect $y, z$ with $|y - z| \leq L/\log L$. By Theorem 5 in Cox and Durrett (2002), if $|y_L - z_L| \geq L/\log L$, then

$$P(|\zeta_{h_L t - s}(\{y_L, z_L\})| = 1) = 1 - \exp(-t + (s/h_L)) + e_L,$$

where $e_L$ is an error term which depends on $L, y_L, z_L, s, t$ and which goes to 0 uniformly for $|y_L - z_L| \geq L/\log L$ and $s \leq t$ in any finite interval. This error term may change from line to line. Using this in the previous equation, we have

$$\int_0^{h_L t} P(\tau = \tau_{kl} \in ds, \tau_{ij} \leq h_L t)$$
$$= \int_0^{h_L t} (1 - \exp(-t + (s/h_L))) P(\tau = \tau_{kl} \in ds) + e_L.$$

Integrating by parts and changing variables, we obtain

$$(4.4) \quad \begin{aligned} \int_0^{h_L t} P(\tau = \tau_{kl} \in ds)\left(1 - \exp\left(-t + \frac{s}{h_L}\right)\right) \\ = \int_0^{h_L t} \frac{1}{h_L} \exp\left(-t + \frac{s}{h_L}\right) P(\tau = \tau_{kl} \leq s) \, ds \\ = \int_0^t e^{-(t-u)} P(\tau = \tau_{kl} \leq h_L u) \, du. \end{aligned}$$

Combining (4.3) and (4.4) yields

$$P(H_t(i,j)) = P(F_t(i,j)) + \sum_{\{k,l\} \neq \{i,j\}} e^{-t} \int_0^t e^s P(F_s(k,l)) \, ds + e_L.$$

Using Theorem 5 in Cox and Durrett (2002) again, $P(H_t(i,j)) \to 1 - e^{-t}$ as $L \to \infty$, which yields

$$(4.5) \quad 1 - e^{-t} = P(F_t(i,j)) + \sum_{\{k,l\} \neq \{i,j\}} e^{-t} \int_0^t e^s P(F_s(k,l)) \, ds + e_L.$$

Summing over all pairs $i, j$ yields

$$\binom{n}{2}(1 - e^{-t}) = q(t) + \left[\binom{n}{2} - 1\right] e^{-t} \int_0^t e^s q(s) \, ds + e_L.$$

It follows [see page 365 of Cox and Griffeath (1986) for details] that $q(t)$ converges to $u(t)$, the solution of

$$\binom{n}{2}(1 - e^{-t}) = u(t) + \left[\binom{n}{2} - 1\right] e^{-t} \int_0^t e^s u(s) \, ds.$$

Rearranging we have

$$e^t u(t) - \binom{n}{2}(e^t - 1) = -\left[\binom{n}{2} - 1\right] \int_0^t e^s u(s) \, ds.$$

Differentiating gives

$$e^t u(t) + e^t u'(t) - \binom{n}{2} e^t = -\left[\binom{n}{2} - 1\right] e^t u(t),$$

which is equivalent to

$$u'(t) = -\binom{n}{2} u(t) + \binom{n}{2}.$$

Since $u(0) = 0$, solving gives $u(t) = 1 - \exp(-\binom{n}{2}t)$. □

While the last calculation is fresh in the reader's mind, we check the claim that when there are $n$ lineages, all $\binom{n}{2}$ coalescences are equally likely. To do this we go back to (4.5). Adding and subtracting $P(F_s(i,j))$ inside the integral,

$$P(F_t(i,j)) - e^{-t} \int_0^t e^s P(F_s(i,j)) \, ds$$

$$= 1 - e^{-t} - e^{-t} \int_0^t e^s q(s) \, ds - e_L.$$

It follows that $P(F_t(i,j))$ converges to $f(t)$, the solution of

$$f(t) - e^{-t} \int_0^t e^s f(s) \, ds = 1 - e^{-t} - e^{-t} \int_0^t e^s u(s) \, ds.$$

Since the limit is independent of $i, j$, it follows that $f(t) = u(t)/\binom{n}{2}$.

PROOF OF THEOREM 2. Lemma 4.2 gives the result for $k = n$ since $P_n(D_t = n) = \exp(-\binom{n}{2}t)$. To prove the result for $k < n$ we use induction on $n$. Theorem 5 of Cox and Durrett (2002) gives the result for $n = 2$. Breaking things down according to the time of the first coalescence, we can write for $B \in \mathcal{G}(L, n, 1)$,

(4.6)
$$\begin{aligned}
P(|\zeta_{h_L t}(B)| = k) \\
= \int_0^{h_L t} P(\tau \in ds, |\zeta_{h_L t}(B)| = k) \\
= \int_0^{h_L t} \sum_{A = \{z_1, \ldots, z_{n-1}\}} P(\tau \in ds, \zeta_s(B) = A) P(|\zeta_{h_L t - s}(A)| = k).
\end{aligned}$$

By Lemma 4.1 it is enough to consider sets $A \in \mathcal{G}(L, n-1, 1)$. The induction hypothesis gives us

$$P(|\zeta_{h_L(t-s)}(A)| = k) = P_{n-1}(D_{t-s} = k) + e_L,$$

where $e_L \to 0$ uniformly for all $A \in \mathcal{G}(L, n-1, 1)$ and $0 \le s \le t$ in any finite interval. Applying the last result again and a change of variables, the quantity on the right-hand side of (4.6) becomes

$$\int_0^t P(\tau/h_L \in ds) P_{n-1}(D_{t-s} = k) + e_L.$$

By Lemma 4.2 we know that

$$P(\tau \le h_L s) = 1 - \exp\left(-\binom{n}{2} s\right) + e_L.$$

Since $s \to P_{n-1}(D_{t-s} = k)$ is continuous, we have

$$P(|\zeta_{h_L t}(B)| = k) \to \int_0^t \binom{n}{2} \exp\left(-\binom{n}{2} s\right) P_{n-1}(D_{t-s} = k) \, ds.$$

The right-hand side is $P_n(D_t = k)$, so the proof of Theorem 2 is complete. □

**5. Proof of Theorem 3.** In view of Theorem 2, it is enough to prove the result for times $L^{2\gamma}/2\nu$ with $0 \le \gamma \le 1$ and show that the ending configuration satisfies the hypotheses of Theorem 2. The second conclusion follows from Lemma 3.7. For the first, it is enough to establish the result up to time $L^2/(2\nu(\log L)^4) = L^{2\kappa_L}/(2\nu)$, where $\kappa_L = 1 - (2 \log \log L)/\log L$, for then Lemma 3.8 implies no collisions occur in $[L^2/(2\nu(\log L)^4), L^2]$.

By rotating the torus we can suppose that $0 \in A$. By the first calculation in the proof of Theorem 1, we can consider the problem on $\mathbb{Z}^2$. So we redefine the following sets as subsets of $\mathbb{Z}^2$. Let

$$\mathcal{G}(L, n, c, \delta) = \{A = \{x_1, \ldots x_n\} : \forall i, x_i \in \mathbb{Z}^2, \forall i \ne j, |x_i - x_j| \in \Gamma(L, c, \delta)\},$$

where $\Gamma(L, c, \delta) = (L^\delta/\log L, c\delta L^\delta \log L)$. Let $\tau$ be the first coalescing time of any two of the $n$ lineages started in $x_1, \ldots, x_n$ and let $\tau_{ij}$ be the coalescing time of the lineages started in $x_i$ and $x_j$. For convenience we let

$$\eta = \frac{\log(2\nu\tau)}{2\log L}, \qquad \eta_{12} = \frac{\log(2\nu\tau_{12})}{2\log L}.$$

LEMMA 5.1. *Let $\zeta_t$ be started with four lineages in $A = \{x_1, \ldots, x_4\}$. Then as $L \to \infty$, uniformly in $A \in \mathcal{G}(L, 4, c, \beta)$,*

$$\int_\beta^{\kappa_L} P(\eta = \eta_{12} \in d\delta, |X_{L^{2\delta}/2\nu}(x_1) - X_{L^{2\delta}/2\nu}(x_3)| \notin \Gamma(L, c+1, \delta)) \to 0,$$

$$\int_\beta^{\kappa_L} P(\eta = \eta_{12} \in d\delta, |X_{L^{2\delta}/2\nu}(x_3) - X_{L^{2\delta}/2\nu}(x_4)| \notin \Gamma(L, c+1, \delta)) \to 0.$$

PROOF. We could repeat the proof of Lemma 1 in Cox and Griffeath (1986), but the following argument is simpler. As in the previous section, we prove only the first statement, since the proof of the second statement is similar. The law of the iterated logarithm implies that

$$P(|X_{t/2\nu}(x_i) - x_i| > \tfrac{1}{4} t^{1/2} \log t \text{ for some } t \geq L^{2\beta}) \to 0.$$

Since $|x_i| \leq c\beta L^\beta \log L \leq \frac{c}{2} t^{1/2} \log t$ for $t \geq L^{2\beta}$, it follows that

$$P\bigg(|X_{t/2\nu}(x_i) - X_{t/2\nu}(x_j)| > \frac{c+1}{2} t^{1/2} \log t \text{ for some } t \geq L^{2\beta}\bigg) \to 0.$$

To show that the particles do not end up too close together, we use the approach of Lemma 4.1. Breaking things down according to the locations of the particles, we want to estimate

$$\int_\beta^{\kappa_L} \sum_{y,z\,:\,|y-z|\leq L^\delta/\log L} P(\eta_{12} \in d\delta, X_{L^{2\delta}/2\nu}(x_1) = y) P(X_{L^{2\delta}/2\nu}(x_3) = z).$$

By the local central limit theorem, the sum over $z$ is at most

$$\bigg(\frac{L^\delta}{\log L}\bigg)^2 \cdot \frac{C}{L^{2\delta}} = \frac{C}{(\log L)^2}.$$

Since

$$\int_\beta^{\kappa_L} \sum_y P(\eta_{12} \in d\delta, X_{L^{2\delta}/2\nu}(x_1) = y) \leq 1,$$

the desired result follows. □

LEMMA 5.2. *As $L \to \infty$,*

$$\sup_{\beta_0 \leq \beta \leq \gamma \leq \kappa_L} \sup_{A \in \mathcal{G}(L,n,c,\beta)} \bigg| P_A(|\zeta_{L^{2\gamma}/2\nu}| = n) - \bigg(\frac{\beta+\alpha}{\gamma+\alpha}\bigg)^{\binom{n}{2}} \bigg| \to 0.$$

PROOF. The proof is a modification of the proof of Proposition 2 of Cox and Griffeath (1986). The case $n = 2$ is covered by Theorem 1. We consider now the case $n > 2$. We need the notation

$$H_\gamma(i,j) = \{\tau_{ij} \leq L^{2\gamma}/2\nu\},$$
$$F_\gamma(i,j) = \{\tau = \tau_{ij} \leq L^{2\gamma}/2\nu\},$$
$$q(\gamma) = P(\tau \leq L^{2\gamma}/2\nu).$$

The estimates in (3.5) and Lemma 3.5 imply that

$$P(\tau \leq L^{2\beta}/2\nu) + P(L^{2\gamma}/4\nu \leq \tau \leq L^{2\gamma}/2\nu) \leq e_L,$$

where here and in what follows $e_L$ is a quantity which depends on $L, A, \beta, \gamma$ and which tends to 0 uniformly for $A \in \mathcal{G}(L, n, c, \beta)$ and $\beta_0 \leq \beta \leq \gamma \leq \kappa_L$. Thus we have

$$
\begin{aligned}
P(H_\gamma(i,j)) = e_L &+ P\left(\tau = \tau_{ij} \leq \frac{L^{2\gamma}}{2\nu}\right) \\
&+ \sum_{\{k,l\} \neq \{i,j\}} \int_{L^{2\beta}/2\nu}^{L^{2\gamma}/4\nu} P\left(\tau = \tau_{kl} \in ds, s < \tau_{ij} \leq \frac{L^{2\gamma}}{2\nu}\right).
\end{aligned}
\tag{5.1}
$$

Letting $\gamma' = \gamma - (\log 2)/(2 \log L)$ so that $L^{2\gamma}/4\nu = L^{2\gamma'}/2\nu$, the $k, l$ term in the last sum is

$$
\int_\beta^{\gamma'} \sum_{y,z} P(\eta = \eta_{kl} \in d\delta, X_{L^{2\delta}/2\nu}(x_i) = y, X_{L^{2\delta}/2\nu}(x_j) = z) \\
\times P(|\zeta_{(L^{2\gamma} - L^{2\delta})/2\nu}(\{y,z\})| = 1).
\tag{5.2}
$$

By Lemma 5.1 we can suppose $|y - z| \in \Gamma(L, c+1, \delta)$. Noting that when $\delta \leq \gamma'$ we have $L^{2\gamma} - L^{2\delta} \geq L^{2\gamma}/2$, and using Theorem 1,

$$
P(|\zeta_{(L^{2\gamma} - L^{2\delta})/2\nu}(\{y,z\})| = 1) = 1 - \frac{\delta + \alpha}{\gamma + \alpha} + e'_L,
$$

where $e'_L \to 0$ uniformly for $|y - z| \in \Gamma(L, c+1, \delta)$ and $\beta \leq \delta \leq \gamma'$. Using this and then replacing the upper limit $\gamma'$ by $\gamma$, we conclude that the quantity in (5.2) is

$$
\int_\beta^\gamma \left(1 - \frac{\delta + \alpha}{\gamma + \alpha}\right) P(\eta = \eta_{kl} \in d\delta) + e_L.
$$

Integrating by parts we obtain

$$
\int_\beta^\gamma \frac{1}{\gamma + \alpha} P\left(\tau = \tau_{kl} \leq \frac{L^{2\delta}}{2\nu}\right) d\delta + e_L = \frac{1}{\gamma + \alpha} \int_\beta^\gamma P(F_\delta(k,l)) \, d\delta + e_L.
$$

Using this in (5.1) yields

$$
P(H_\gamma(i,j)) = P(F_\gamma(i,j)) + \sum_{\{k,l\} \neq \{i,j\}} \frac{1}{\gamma + \alpha} \int_\beta^\gamma P(F_\delta(k,l)) \, ds + e_L.
$$

Since $|x_i - x_j| \in \Gamma(L, c, \beta)$, Theorem 1 implies

$$
1 - \frac{\beta + \alpha}{\gamma + \alpha} = P(F_\gamma(i,j)) + \sum_{\{k,l\} \neq \{i,j\}} \frac{1}{\gamma + \alpha} \int_\beta^\gamma P(F_\delta(k,l)) \, d\delta + e_L.
\tag{5.3}
$$

Summing over all pairs $i, j$,

$$
\binom{n}{2}\left(1 - \frac{\beta + \alpha}{\gamma + \alpha}\right) = q(\gamma) + \left[\binom{n}{2} - 1\right] \frac{1}{\gamma + \alpha} \int_\beta^\gamma q(\delta) \, d\delta + e_L.
$$

It follows that $q(t)$ converges to $u(t)$, the solution of

$$\binom{n}{2}\left(1 - \frac{\beta + \alpha}{\gamma + \alpha}\right) = u(\gamma) + \left[\binom{n}{2} - 1\right]\frac{1}{\gamma + \alpha}\int_\beta^\gamma u(\delta)\,d\delta.$$

This leads to

$$u(\gamma) = 1 - \left(\frac{\beta + \alpha}{\gamma + \alpha}\right)^{\binom{n}{2}}.$$

From this it follows [see page 365 of Cox and Griffeath (1986) for more details] that

$$q(\gamma) = 1 - \left(\frac{\beta + \alpha}{\gamma + \alpha}\right)^{\binom{n}{2}} + e_L.$$

This completes the proof of Lemma 5.2. □

Again we pause to check the claim that when there are $n$ lineages, all $\binom{n}{2}$ coalescences are equally likely. We proceed in the same way as in the argument after the proof of Lemma 4.1. We go back to (5.3) and add and subtract $P(F_\delta(i,j))$ inside the integral. It follows that $P(F_\gamma(i,j))$ converges to $f(\gamma)$, the solution of

$$f(\gamma) - \frac{1}{\gamma + \alpha}\int_\beta^\gamma f(\delta)\,d\delta = 1 - \frac{\beta + \alpha}{\gamma + \alpha} - \frac{1}{\gamma + \alpha}\int_\beta^\gamma q(\delta)\,d\delta.$$

Since the limit is independent of $i, j$, it follows that $f(\gamma) = u(\gamma)/\binom{n}{2}$.

PROOF OF THEOREM 3. Lemma 5.2 gives the result for $k = n$. To prove the result for $k < n$, we use induction on $n$. Theorem 1 gives the result for $n = 2$ since

$$P_n(D_{\log((\gamma+\alpha)/(\beta+\alpha))} = n) = \left(\frac{\beta + \alpha}{\gamma + \alpha}\right)^{\binom{n}{2}}.$$

As before,

$$P(\tau \leq L^{2\beta}/2\nu) + P(L^{2\gamma}/4\nu \leq \tau \leq L^{2\gamma}/2\nu) \leq e_L$$

is a quantity that tends to 0 uniformly for $B \in \mathcal{G}(L, n, c, \beta)$ and $\beta_0 \leq \beta \leq \gamma \leq \kappa_L$. So letting $\gamma' = \gamma - (\log 2)/(2\log L)$ as before, we have

$$(5.4) \qquad P(|\zeta_{L^{2\gamma}/2\nu}(B)| = k) = e_L + \int_\beta^{\gamma'} P(\eta \in d\delta, |\zeta_{L^{2\gamma}/2\nu}(B)| = k).$$

As in the proof of Theorem 2, we can write the integral as

$$\int_\beta^{\gamma'} \sum_{A = \{z_1, \ldots z_{n-1}\}} P(\eta \in d\delta, \zeta_{L^{2\delta}/2\nu}(B) = A)P(|\zeta_{(L^{2\gamma} - L^{2\delta})/2\nu}(A)| = k).$$

STEPPING STONE MODEL. II 29

By Lemma 5.1 it is enough to consider sets $A \in \mathcal{G}(L, n-1, c+1, \delta)$, for which we know by the induction hypothesis that

$$P(|\zeta_{(L^{2\gamma}-L^{2\delta})/2\nu}(A)| = k) = P_{n-1}(D_{\log((\gamma+\alpha)/(\delta+\alpha))} = k) + e'_L,$$

where $e'_L \to 0$ uniformly for all $A \in \mathcal{G}(L, n-1, c+1, \delta)$ and $\beta \leq \delta \leq \gamma' \leq \gamma \leq \kappa_L$. By Lemma 5.2 we know that

$$P\left(\tau \leq \frac{L^{2\delta}}{2\nu}\right) = 1 - \left(\frac{\beta+\alpha}{\delta+\alpha}\right)^{\binom{n}{2}} + e_L.$$

Since $\delta \to P_{n-1}(D_{\log((\gamma+\alpha)/(\delta+\alpha))} = k)$ is continuous, we obtain

$$P(|\zeta_{L^{2\gamma}/2\nu}(B)| = k) \to \int_\beta^\gamma \frac{\binom{n}{2}(\beta+\alpha)^{\binom{n}{2}}}{(\delta+\alpha)^{\binom{n}{2}+1}} P_{n-1}(D_{\log((\gamma+\alpha)/(\delta+\alpha))} = k) \, d\delta.$$

Changing variables $\delta = (\beta+\alpha)e^s - \alpha$, $d\delta = (\beta+\alpha)e^s \, ds$, we see that the above

$$\int_0^{\log((\gamma+\alpha)/(\beta+\alpha))} \binom{n}{2} e^{-s\binom{n}{2}} P_{n-1}(D_{\log((\gamma+\alpha)/(\beta+\alpha))-s} = k) \, ds$$
$$= P_n(D_{\log((\gamma+\alpha)/(\beta+\alpha))} = k),$$

which completes the proof of Theorem 3. $\square$

I. Zähle
Mathematische Institut
University of Erlangen-Nuernberg
Bismarckstrasse 1 1/2
910054 Erlangen
Germany
e-mail: zaehle@ami.uni-erlangen.de
url: http://www.mi.uni-erlangen.de/~zaehle/index_e.html

J. T. Cox
Department of Mathematics
Syracuse University
Syracuse, New York 13210
USA
e-mail: jtcox@mailbox.syr.edu
url: http://web.syr.edu/~jtcox/index.html

R. Durrett
Department of Mathematics
Cornell University
523 Malott Hall
Ithaca, New York 14853
USA
e-mail: rtd1@cornell.edu
url: www.math.cornell.edu/~durrett